\documentclass[12pt]{article}

\usepackage{amsmath,amsthm}
\usepackage{amssymb}
\usepackage{color}
\usepackage{xspace}
\usepackage[colorlinks=true,
linkcolor=green,
filecolor=brown,
citecolor=green]{hyperref}

\def\si{\sigma}
\def\p{\ensuremath{\mathcal P}\xspace}

\def\mbf#1{\mathchoice{\hbox{\boldmath $\displaystyle #1$}}
        {\hbox{\boldmath $\textstyle #1$}}
        {\hbox{\boldmath $\scriptstyle #1$}}
        {\hbox{\boldmath $\scriptscriptstyle #1$}}} 

\begin{document}
\newtheorem{theorem}{Theorem}
\newtheorem{defn}[theorem]{Definition}
\newtheorem{lemma}[theorem]{Lemma}
\newtheorem{prop}[theorem]{Proposition}
\newtheorem{cor}[theorem]{Corollary}
\begin{center}
{\Large
An involution on set partitions         \\[2mm] 
}

David Callan  \\

April 5, 2022
\end{center}

\begin{abstract}
We present an involution on set partitions that 
interchanges two statistics related to relative size of block entries and use it to establish an equidistribution on objects counted by the Bessel numbers.

\end{abstract}

\section{Introduction} \label{intro}
Claesson \cite{claesson01} showed bijectively that permutations of $[n]$ avoiding both of the vincular patterns $\overline{12}3$ and $1\overline{23}$ are equinumerous with nonoverlapping set partitions of $[n]$, counted by the Bessel numbers, sequence \htmladdnormallink{A006789}{http://oeis.org/A006789} in OEIS \cite{oeis}.
Mansour and Shattuck recently determined the distribution of the ``last entry'' statistic on these $\{\overline{12}3,\,1\overline{23}\}$-avoiders \cite{li} in the course of enumerating $\overline{23}41$-avoiders. 
Here are the first few values of the number $v_{n,k}$ of $\{\overline{12}3,\,1\overline{23}\}$-avoiders on $[n]$ with last entry = $k$:
\[
\begin{array}{c|ccccccc}
n^{\textstyle{\,\backslash \,k}} & 1 & 2 & 3 & 4  & 5  & 6 &7 \\
\hline 
	1&    1 &   & & & & &  \\
 	2&    1 & 1 & & &  & &  \\
	3&    2 & 2 & 1 & & & &  \\
	4&    5 & 5  & 3 & 1 & &  &  \\ 
	5&    14 & 14 & 9 & 5 & 1 &  &  \\ 
	6&    43 & 43 & 29 & 18 & 9 & 1  &  \\ 
	7&   143 & 143 & 100 & 66 & 39 & 17 & 1  
 \end{array}
\]\\[-4mm]
\centerline{Table of values of $v_{n,k}$}\\[4mm]
and $v_{n,k}$ is defined recursively  \cite{li} by
\begin{eqnarray*}
v_{n,n}\  = & 1    & \qquad \textrm{ for $ n \ge 1$}   \\
v_{n,1}\  = & \sum_{i=1}^{n-1}v_{n - 1, i}  & \qquad \textrm{ for $ n \ge 2$}   \\ 
v_{n,k}\  = & \sum_{i=k}^{n-1}v_{n-1,i} + 
  \sum_{i=k+1}^n \sum_{d=2}^{k} \binom{k - 2}{d - 2} v_{n - d, i - d}  & \qquad \textrm{ for $ 2\le k \le n - 1 $}\, .
\end{eqnarray*}

Under Claesson's bijection, the previously mentioned last entry translates to a statistic $Y$ on nonoverlapping partitions defined in the following section, that turns out to be a variant of the ``minimax'' statistic, $X$, defined as the minimum taken over the maximum entry in each block.
We wish to show that 
$X$ and $Y$ have the same distribution, in fact, a symmetric joint distribution. For this purpose, after some preliminaries in Section \ref{prelim}, Section \ref{invol} presents a size-preserving involution 
$\si$ on set partitions that interchanges the values of $X$ and $Y$ and preserves nonoverlapping partitions.

\section{Preliminaries} \label{prelim}

It is convenient for us to aways write set partitions in a standard form, each block decreasing and blocks arranged in order of increasing first entry, as in  31/62/7/854 with 4 blocks. With this form, the ``minimax'' statistic $X$ is the first entry of the first block.

The \emph{span} of a block $B$ in a set partition $P$ of $[n]$ is the smallest interval $\overline{B}$ of integers containing $B$. Thus $\overline{\{854\} } = [4,8]:=\{4,5,6,7,8\}$. A partition $P$ is \emph{nonoverlapping} if for each pair of blocks $B_1$ and $B_2$ in $P$, the spans  $\overline{B_1}$ and  
$\overline{B_2}$ are either disjoint or one contains the other. Thus 2/43/651/87 is nonoverlapping but 31/62/7/854 is not. Note that if $B_1$ is a singleton block, then $B_1=\overline{B_1}$ is necessarily either disjoint from, or contained in, $\overline{B_2}$ and so only non-singleton blocks need be checked for the nonoverlapping property.

The statistic $Y$ mentioned in the Introduction is defined on all set partitions $P$ as follows:
$Y(P)=  1$ if $\{1\}$ is a singleton block in $P$, and otherwise  
$Y(P)= \min\{r,s\}$   where $r$ is the minimum taken over the maximum entry in each non-singleton block of $P$ and $s$ is the second smallest entry in the block containing 1.
Thus, $Y(1/32)=1$ and $Y(3/4/652/7/981)=\min\{r\!=\!\min\{6,9\},s\!=\!8\}=6$.

\section{The involution $\mbf{\sigma}$ on set partitions} \label{invol}
Let $\p_n$ denote the set of partitions of $[n],\ n\ge 1$, written in standard form.
For $P\in \p_n$, if  $X(P)=Y(P)$, set $\si(P)=P$. If  $X(P)<Y(P)$, let $r$ denote the first entry of the first nonsingleton block in $P$, and let $s$ denote the (left) neighbor of 1, so that $Y(P)=\min(r,s)$.
\begin{itemize}
\vspace*{-2mm}
\item  If $r>s$ as in, for example, $P=3/4/7/852/961$ with $X=3,\ r=8$ and $Y=s=6$, delete all the 
initial singleton blocks whose (sole) entry is $<s$ and reinsert their 
entries, in decreasing order, between $s$ and 1 in the block containing 1. 
Then transfer $s$ from its block to form a new singleton block at the far left: $\sigma(P) = 6/7/852/9431$ in the example. 
A similar example, where $r$ and $s$ are in the same block, is $\sigma(2/431) = 3/421$.  

In this case, $\sigma(P)$ will begin with a singleton block.
\vspace*{-1mm}
\item  If $r\le s$, as in $P=3/4/652/7/981$ with $X=3,\ s=8$ and $Y=r=6$, delete all the 
initial singleton blocks and reinsert their 
entries, necessarily $<r$, between $s$ and 1 in the block containing 1: 
here $\sigma(P) = 652/7/98431$. A similar example, where $r=s$, is $\sigma(2/3/4/51) = 54321$.  

In this case, $\sigma(P)$ will begin with a nonsingleton block.
\end{itemize}

Clearly, with this definition, $\si$ interchanges the values of $X$ and 
$Y$ and so sends
$\{P\in \p_n:X(P)<Y(P)\}$ to $\{P\in \p_n:X(P)>Y(P)\}$ for $n\ge 3$.
It is easy to see this latter map is onto and invertible and so (along with $\si(P)=P$ when $X(P)=Y(P)$\,) 
defines an involution $\si$ on $\p_n$ that interchanges the values of $X$ and $Y$.  We leave the details to the reader. The map $\si$ preserves the span of the nonsingleton blocks in each partition $P\in \p_n$ and so $\si$ preserves the set of nonoverlapping partitions of $[n],\ n\ge 1$.

Department of Statistics, University of Wisconsin-Madison

\end{document}